\documentclass[11pt]{amsart}

\usepackage[margin=1.15in]{geometry}

\usepackage{amsmath}
\usepackage{amssymb}
\usepackage{amsthm}
\usepackage{enumerate}
\usepackage{mathdots}

\usepackage[active]{srcltx}

\usepackage{amsfonts}
\usepackage{enumerate}
\usepackage{bbm}

\usepackage[dvips]{color}

\definecolor{refkey}{gray}{0.5}
\definecolor{labelkey}{gray}{0.5}

\newcommand{\nC}{\mathbb C}

\newcommand{\nN}{\mathbb N}

\newcommand{\M}[1]{{#1}}

\newtheorem*{theorem*}{Theorem}

\theoremstyle{definition}

\begin{document}

\title[The main theorem on Bezoutians]{A proof of the main theorem on Bezoutians}

\author{Branko \'{C}urgus}
\author{Aad Dijksma}

\date{\today}

\begin{abstract}
We give a self-contained proof that the nullity of the Bezoutian matrix associated with a pair of polynomials $f$ and $g$ equals the number of their common zeros counting multiplicities.
\end{abstract}

\maketitle

With two polynomials $f$ and $g$ and $n = \max\{\deg f, \deg g\}$ we associate an $n \times n$ matrix $B$,
called the Bezoutian, and a $2n \times 2n$ matrix $R$, called the resultant. Their defining relations are given by \eqref{eqdB} and \eqref{eqRV} below, respectively. In terms of the coefficients of $f$ and $g$ they are given by \eqref{eqBHT} and \eqref{eqdR}. In this note we give a simple and self-contained proof of the equalities
\begin{equation}\label{desired=}
\dim \ker B = \dim \ker R = \deg \gcd (f,g),
\end{equation}
where $\gcd$ stands for greatest common divisor. H.~Wimmer in \cite{Wimmer} attributes this result to Jacobi who in 1836 showed that the singularity of what we call the Bezoutian implies the existence of a common factor of $f$ and $g$.  More contemporary proofs of \eqref{desired=} can be found in the recent books \cite[Theorems~21.10 and~21.11]{Dym} by H.~Dym and \cite[Theorem~8.30]{Fuhrmann} by P.~Fuhrmann.  In the Introduction to \cite[Chapter~21]{Dym} it is shown that $\dim \ker B \geq \deg\gcd(f,g)$ by using the defining formula for $B$, differentiation and chains of vectors. That equality prevails is then proved by using these chains and the so-called Barnett identity: $B=H_fg(C_f)$, where $H_f$ is the Hankel matrix for $f$ defined below and $C_f$ is the companion matrix of $f$.  In \cite{Fuhrmann} the matrix $B$ is expressed in terms of a matrix representation of $g(S_f)$, where $S_f$ is the shift operator in the space $X_f$ of polynomials modulo $f$,  relative to two suitably chosen bases in $X_f$. In view of
\cite[Corollary~8.29]{Fuhrmann} this formula is closely related to the Barnett identity. In this note we do not resort to this identity. Our approach, we think, is more direct. Of course some of the formulas derived below also appear in \cite[Chapter 21]{Dym} and \cite[Chapter~8]{Fuhrmann}.
Our proofs of these formulas are different. For a survey of results related to Bezoutians see \cite[Fact~4.8.6]{Bernstein} in the encyclopedic book by D.S.~Bernstein and for
applications of Bezoutians in  numerical linear algebra and system theory,
see for example \cite{Fiedler} and \cite{Hinrichsen}, respectively.

\medskip

\noindent {\bf 1.~Notation and basic notions.} The vector space of all polynomials with coefficients in $\nC$ and in the variable $z$ is denoted by $\nC[z]$. Its  Cartesian square is denoted by $\nC^2[z]$. For $n\in \nN$,  $\nC[z]_{<n}$ denotes the subspace of $\nC[z]$ of all polynomials of degree strictly less than $n$. This space has dimension $n$. Similarly, $\nC^2[z]_{<n}$ denotes the Cartesian square of $\nC[z]_{<n}$.

We use $I$ to denote the identity matrix, $Z$  the reverse identity and $N$ the nilpotent Jordan block:
\begin{equation*} 
{Z}: = \begin{bmatrix}
0 &  \cdots  & 1 \\[-2pt]
\vdots &  \iddots  & \vdots \\
1 &  \cdots  & 0
\end{bmatrix}, \qquad \qquad
\addtolength{\arraycolsep}{-4pt}
{N}: = \begin{bmatrix}
0 & 1 &  \cdots  & 0 \\[-5pt]
\vdots &  \ddots &  \ddots  & \vdots \\[-8pt]
0 &  & \ddots  & 1 \\[-2pt]
0 & 0 &  \cdots  & 0
\end{bmatrix}. \addtolength{\arraycolsep}{0mm}
 \end{equation*}

For a polynomial
\[
f(z) = f_0 + f_1 z + \cdots f_n z^n
\]
in $\nC[z]$ we define two $n\times n$ matrices, one Hankel and one Toeplitz,  associated with $f$ as follows:
\begin{alignat*}{2} 
{H}_{f} &: = \begin{bmatrix}
f_1 &  \cdots  & f_n \\
\vdots & \iddots  & \vdots \\
f_n & \cdots & 0
\end{bmatrix}, &
\qquad
{T}_{f} &: =
\begin{bmatrix}
f_0  & \cdots & f_{n-1} \\
\vdots & \ddots  & \vdots  \\
0 &  \cdots  & f_0
\end{bmatrix}. \\
\intertext{Since the left-multiplication by $Z$ reverses the rows, it turns a Hankel matrix into a Toeplitz and vice versa:}
{Z} {H}_{f} &: = \begin{bmatrix}
f_n & \cdots  & 0 \\
\vdots &  \ddots  & \vdots \\
f_1 & \cdots  & f_n \\
\end{bmatrix}, &
\qquad
{Z} {T}_{f} &: =
\begin{bmatrix}
0 & \cdots  & f_0\\
\vdots &  \iddots  & \vdots  \\
f_0  & \cdots & f_{n-1}
\end{bmatrix}.
 \end{alignat*}
As each Hankel matrix is symmetric, we have ${Z} {T}_{f} = \bigl({Z} {T}_{f}\bigr)^\top = T_f^\top Z$, where the superscript $\phantom{}^\top$ is used to denote a matrix transpose. Consequently,
\begin{equation}\label{eqTT}
 {T}_{f}^\top =  ZT_fZ.
\end{equation}

The vector space of upper (lower) triangular Toeplitz matrices is spanned by the identity $I$ and the powers of $N$ ($N^\top$, respectively). Therefore, the upper (lower) triangular Toeplitz matrices form a commutative algebra. In particular for polynomials $f$ as above and $g(z)=g_0 + g_1 z + \cdots g_n z^n$  we have
\begin{equation}\label{eqcom1}
T_{f}T_{g} = T_{g}T_{f}, \qquad H_{f} Z H_{g}  = H_{g} Z H_{f},
\end{equation}
where the last equality follows from $(Z H_{f}) (Z H_{g}) = (Z H_{g}) (Z H_{f})$.

For $n \in \nN$ and $z \in \nC$  we denote by $V_n(z)$ the $n\times 1$ column vector
\[
V_n(z) = \begin{bmatrix}
1 & z & \cdots & z^{n-1}
\end{bmatrix}^\top.
\]
This notation is convenient as it provides a compact way of writing polynomials. For example, a polynomial $a(z,w)$ in two variables $z$ and $w$ can be written as:
\[
a(z,w) = \sum_{j,k=0}^{n-1} a_{jk} z^j  w^k = {V}_{n}(z)^\top {A} {V}_{n}(w),
\]
where ${A}$ is the $n\times n$ coefficient matrix $\bigl[a_{jk}\bigr]_{j,k=0}^{n-1}$ of $a(z,w)$.

The {\em resultant} $R$ of the polynomials $f$ and $g$ is the $2n\times 2n$ matrix given as a $2\times 2$ block matrix:
\begin{equation}\label{eqdR}
\addtolength{\arraycolsep}{-2pt}
R = \begin{bmatrix}
{T}_{f} & Z  H_{f}  \\[8pt]
{T}_{g} & Z  H_{g}
\end{bmatrix}.
\addtolength{\arraycolsep}{0pt}
\end{equation}
Notice that the action of $R$ on $V_{2n}(z)$ is particularly simple:
\begin{equation}\label{eqRV}
RV_{2n}(z) =\begin{bmatrix}
f(z) V_{n}(z) & g(z) V_{n}(z)
\end{bmatrix}^\top.
\end{equation}

Next we define the Bezoutian $B$ of $f$ and $g$. First consider the polynomial $f(z)g(w) - f(w)g(z)$ in two variables. Since this polynomial vanishes for all $w = z \in \nC$, there exists a polynomial $b(z,w)$ in two variables such that
\[
f(z)g(w) - f(w)g(z) = (z-w)b(z,w) \quad\text{for all} \quad z,w \in \nC.
\]
The {\em Bezoutian} $B $ of $f$ and $g$ is the $n\times n$ coefficient matrix of $b(z,w)$:
\begin{equation}\label{eqdB}
b(z,w) = {V}_{n}(z)^\top {B } {V}_{n}(w), \quad z,w \in \nC.
\end{equation}

The null space or kernel of a matrix (or a linear transformation) $A$ is denoted by $\ker A$. Its dimension is called the {\em nullity} of $A$.

\medskip
\noindent {\bf 2.~A connection between $R$ and $B$.} To establish a connection between $R$ and $B$ we consider the polynomial $(z^n-w^n)b(z,w)$ and we find two ways of representing its coefficient matrix. To find the first representation we use the standard identity
\[
z^n - w^n = (z-w)\sum_{j=0}^{n-1} z^{n-1-j}w^j = (z-w) {V}_n(z)^\top {Z}  {V}_n(w), \quad z,w \in \nC,
\]
matrix algebra and \eqref{eqRV}:
\allowdisplaybreaks{%
\begin{align*}
(z^n-w^n)b(z,w)
& = (z-w)b(z,w) \M{V}_n(z)^\top\M{Z}\M{V}_n(w) \\
& =
\bigl(f(z)g(w)- g(z)f(w) \bigr)\M{V}_n(z)^\top\M{Z}\M{V}_n(w) \\
& =
\begin{bmatrix}
f(z) {V}_n(z) \\[8pt]
g(z) {V}_n(z)
\end{bmatrix}^\top
\begin{bmatrix}
g(w) {Z}{V}_n(w) \\[8pt]
-f(w) {Z}{V}_n(w)
\end{bmatrix} \\
& =
\begin{bmatrix}
f(z) {V}_n(z) \\[8pt]
g(z) {V}_n(z)
\end{bmatrix}^\top \begin{bmatrix}
0 & Z \\[8pt]
-Z & 0
\end{bmatrix}\begin{bmatrix}
f(w) {V}_n(w) \\[8pt]
g(w) {V}_n(w)
\end{bmatrix} \\
& = {V}_{2n}(z)^{\!\top}\! {R}^\top\!\! \begin{bmatrix}
0 & Z \\[8pt]
-Z & 0
\end{bmatrix}\! R {V}_{2n}(w).
\end{align*}
}
\!\!The second representation involves the Bezoutian:
\allowdisplaybreaks{%
\begin{align*}
 (z^n-w^n)b(z,w)
& = (z^n-w^n) {V}_n(z)^\top B {V}_n(w)\\
& = \bigl(z^n {V}_n(z)\bigr)^{\!\top}\! B {V}_n(w)\! -\! {V}_n(z)^{\!\top}\! B \bigl( w^n {V}_n(w)\bigr)\\
& =  {V}_{2n}(z)^\top \begin{bmatrix}
0 & 0 \\[8pt]
B & 0 \end{bmatrix} {V}_{2n}(w)
 + {V}_{2n}(z)^\top
 \addtolength{\arraycolsep}{-2pt}
 \begin{bmatrix}
0 & -B \\[8pt]
0 & 0
\end{bmatrix} {V}_{2n}(w)\\
& = {V}_{2n}(z)^\top
\addtolength{\arraycolsep}{-2pt}
\begin{bmatrix}
0 & -B \\[8pt]
B & 0
\end{bmatrix} {V}_{2n}(w).
\addtolength{\arraycolsep}{0pt}
\end{align*}
}
\!\!These two representations of the coefficient matrix of $(z^n-w^n)b(z,w)$ provide a connection between $R$ and $B$:
\begin{equation} \label{eqRandB}
R^\top \begin{bmatrix}
0 & Z \\[8pt]
-Z & 0
\end{bmatrix} R =
\addtolength{\arraycolsep}{-2pt}
\begin{bmatrix}
0 & -B \\[8pt]
 B & 0
\end{bmatrix}.
\end{equation}

On the other hand, using the definition of $R$, \eqref{eqTT} and  \eqref{eqcom1} we obtain
\allowdisplaybreaks{%
\begin{align*}
R^\top \begin{bmatrix}
0 & Z \\[8pt]
-Z & 0
\end{bmatrix} R & =  \begin{bmatrix}
ZT_fZ & ZT_gZ  \\[8pt]
H_f Z & H_g Z
\end{bmatrix}
\begin{bmatrix}
Z T_g & H_g  \\[8pt]
-ZT_f & -H_f
\end{bmatrix}\\
& = \begin{bmatrix}
ZT_fT_g - ZT_gT_f & Z T_f Z H_g - Z T_g Z H_f \\[8pt]
H_f T_g - H_g T_f & H_f Z H_g - H_gZH_f
\end{bmatrix} \\
&= \begin{bmatrix}
0 & -\bigl(H_f T_g - H_g T_f\bigr)^\top \\[8pt]
 H_f T_g - H_g T_f & 0
\end{bmatrix}.
\end{align*}
}
\!\!Together with \eqref{eqRandB}, the last equality yields
\[
\begin{bmatrix}
0 & -\bigl(H_f T_g - H_g T_f\bigr)^\top \\[8pt]
 H_f T_g - H_g T_f & 0
\end{bmatrix} =
\begin{bmatrix}
0 & -B \\[8pt]
 B & 0
\end{bmatrix},
\]
and thus
\begin{equation} \label{eqBHT}
B = H_f T_g - H_g T_f = B^\top.
\end{equation}

\medskip
\noindent {\bf 3.~$R$ and $B $ have the same nullity.}
Equation \eqref{eqRandB} indicates that there is a connection between $\ker R$ and $\ker B$. An even more direct connection between $\ker R$ and $\ker B$ is obtained from \eqref{eqdR}, \eqref{eqBHT} and \eqref{eqcom1} (listed in the order in which they are used) as follows:
\begin{equation} \label{eqaad}
\begin{split}
\begin{bmatrix} I & 0 \\ T_f & ZH_f \end{bmatrix}
 R
 & = \begin{bmatrix} T_f & ZH_f \\ T_f^2+ ZH_f T_g & T_fZH_f + ZH_f ZH_g \end{bmatrix} \\
 & = \begin{bmatrix} T_f & ZH_f \\ ZB + \bigl(T_f+ ZH_g\bigr) T_f & \bigl(T_f + ZH_g\bigr) ZH_f \end{bmatrix} \\
 & = \begin{bmatrix} 0 & I \\ ZB & T_f+ZH_g \end{bmatrix}
\begin{bmatrix} I & 0 \\ T_f & ZH_f \end{bmatrix} \\
& = \begin{bmatrix} 0 & I \\ Z & T_f+ZH_g \end{bmatrix}
\begin{bmatrix} B&0 \\ 0 & I\end{bmatrix}
\begin{bmatrix} I & 0 \\ T_f & ZH_f \end{bmatrix}.
\end{split}
\end{equation}
If we assume that $n = \deg f$, then $H_f$ is invertible, yielding that the first (as well as the last) block matrix in \eqref{eqaad} is invertible. Since the block matrix in \eqref{eqaad} whose antidiagonal  entries are $I$ and $Z$ is also invertible, \eqref{eqaad} implies that $R$ and $B$ have the same nullities:
\begin{equation} \label{eqnulRB}
\dim \ker R = \dim \ker B.
\end{equation}

\medskip
\noindent {\bf 4.~The nullity of $B$ in terms of $f$ and $g$.}
Consider the multiplication operator
$$M: \nC^2[z]_{<n} \to \nC[z]_{<2n}$$
defined by
\[
M \begin{bmatrix}
u \\[2pt] v
\end{bmatrix} = fu+gv, \quad u, v \in \nC[z]_{<n}.
\]
For a characterization of the null space $\ker M$ of $M$  in terms of $f$ and $g$ we need the greatest common divisor $h$ of $f$ and $g$, its degree  $k=\deg h$ and factorizations $f = \hat{f} h$, $g = \hat{g} h$. Then
\begin{equation} \label{eqnulM}
\ker M = \left\{ \begin{bmatrix}
u \\ v \end{bmatrix} \in \nC^2[z] \, : \, u = -\hat{g} q, \, v = \hat{f} q, \, q \in \nC[z]_{<{k}} \right\}.
\end{equation}
The inclusion $\supseteq$ in \eqref{eqnulM} is clear. To prove $\subseteq$, let $u,v \in \nC[z]_{<n}$ and $\begin{bmatrix}
u & v \end{bmatrix}^\top \in \ker M$. Then $fu+gv=0$, implies  $\hat{f}u = -\hat{g}v$. Since $\hat{f}$ and $\hat{g}$ have no common zeros, the last identity yields that there exist polynomials $p$ and $q$ such that $u = \hat{g}p$ and $v=\hat{f}q$. Substituting back to $\hat{f}u = -\hat{g}v$, we get $\hat{f}\hat{g}p = -\hat{g}\hat{f}q$. Hence $p=-q$. Since $\deg v <n$ and $\deg \hat{f} = n - k$, $v = \hat{f} q$ implies $\deg q < k$. This proves \eqref{eqnulM}.

The standard basis for $\nC^2[z]_{<n}$ is
\[
\begin{bmatrix}
1 \\[2pt] 0
\end{bmatrix}, \ldots,
\begin{bmatrix}
z^{n-1} \\[2pt] 0
\end{bmatrix}, \begin{bmatrix}
0 \\[2pt] 1
\end{bmatrix}, \ldots,
\begin{bmatrix}
0 \\[2pt] z^{n-1}
\end{bmatrix},
\]
while the standard basis for $\nC[z]_{<2n}$ is
\[
1, z, \ldots, z^{n-1},z^n, \ldots, z^{2n-1}.
\]
The matrix representation for $M$ with respect to these standard bases is $R^\top$, see \eqref{eqRV}. Therefore the nullity of $R^\top$ is $\dim \ker M$. Since \eqref{eqnulM} yields $\dim \ker M = k$ and the nullity of $R^\top$ equals the nullity of $R$, we have proved that the nullity of $R$ is $k$ and, by \eqref{eqnulRB},
\begin{align*}
\dim \ker B = \dim \ker R = k = \deg h.
\end{align*}

\noindent{\bf 5.~Final remarks.} It was remarked in \cite[p.~318]{Hinrichsen} that the Bezoutian of a pair of polynomials is  defined whenever $n \geq \max\{\deg f,\deg g\}$. We add to this that the same is true for the resultant and that if $n \geq m:=\max\{\deg f,\deg g\}$, then formula \eqref{desired=} has to be replaced by the formula
\begin{equation}\label{general}
\dim \ker B_n = \dim \ker R_{2n} =n-m +\deg \gcd (f,g),
\end{equation}
where, for example, the index $n$ in $B_n$ indicates that $B_n$ has size $n \times n$. Indeed, \eqref{general} follows from \eqref{desired=} and from the equalities
\[
\dim \ker B_n=n-m+ \dim \ker B_m  \quad
\text{and} \quad
\dim \ker R_{2n}=n-m+ \dim \ker R_{2m}.
\]
The first of these two equalities holds because of \eqref{eqdB}, which implies $B_n=\begin{bmatrix} B_m \!\!\! & 0 \\ 0 & 0 \end{bmatrix},$
and the second follows from the reasoning in Section~4 with $k$ in \eqref{eqnulM} replaced by $n-m+k$.

Finally we note that \eqref{general} can be expressed as
\begin{equation}\nonumber
\dim \ker B_n = \dim \ker R_n = \deg \gcd (\bar{f},\bar{g}),
\end{equation}
where
\[
\bar{f}(y,z) = f_0y^n + f_1 y^{n-1} z + \cdots +f_nz^n \ \   \text{and} \ \ \bar{g}(y,z) = g_0y^n + g_1 y^{n-1} z + \cdots +g_nz^n
\]
are homogenizations (in the sense of \cite[page 6-7]{Marshall}) of $f$ and $g$, respectively. If $n > m = \max \{\deg f, \deg g\}$, then $y=0$ is a common zero of $\bar{f}$ and $\bar{g}$ of multiplicity $n-m$. Since the zero $y=0$ of the homogenization is commonly viewed as a ``zero at infinity'' of the original polynomial (see for example \cite[4.4.3]{Blekherman}) we can rephrase our abstract to cover the ``generalized'' Bezoutian $B_n$:  The nullity of the Bezoutian matrix $B_n$ associated with a pair of polynomials $f$ and $g$ equals the number of their common zeros including the ``zero at infinity'' and counting multiplicities.


\begin{thebibliography}{10}

\bibitem{Bernstein}
D.S. Bernstein, Matrix Mathematics: Theory, Facts, and Formulas. Second Edition, Princeton University Press, 2009.

\bibitem{Blekherman}
G.~Blekherman, Nonnegative polynomials and sums of squares, in Semidefinite Optimization and Convex Algebraic Geometry edited by G.~Blekherman, P.A.~Parrilo, R.~Thomas, forthcoming book in the MOS-SIAM Series on Optimization.


\bibitem{Dym}
H.~Dym, Linear algebra in action. Graduate Studies in Mathematics, 78. American Mathematical Society, Providence, RI, 2007.

\bibitem{Fiedler}
M. Fiedler, Special Matrices and Their Applications in Numerical Mathematics. Second Edition, Dover Books on Mathematics, Mineola, NY,  2008.

\bibitem{Fuhrmann}
P.~Fuhrmann, A polynomial approach to linear algebra. Second edition. Universitext. Springer, New York, 2012.

\bibitem{Hinrichsen}
D. Hinrichsen, A. Pritchard, Mathematical systems theory I.
Modelling, state space analysis, stability and robustness.  Texts in Applied Mathematics, 48. Springer, Heidelberg, 2010.

\bibitem{Marshall}
M.~Marshall, Positive polynomials and sums of squares. Mathematical Surveys and Monographs, 146. American Mathematical Society, Providence, RI, 2008.

\bibitem{Wimmer}
H.K. Wimmer, On the history of the Bezoutian and the resultant matrix. Linear Algebra Appl. 128 (1990), 27–-34.


\end{thebibliography}
\end{document}